\pdfoutput=1 
\documentclass[fleqn,12pt]{SelfArx} 

\usepackage[english]{babel} 

\usepackage{lipsum} 

\definecolor{color1}{RGB}{0,0,90} 
\definecolor{color2}{RGB}{0,20,20} 

\usepackage{hyperref} 

\hypersetup{
	hidelinks,
	colorlinks,
	breaklinks=true,
	urlcolor=color2,
	citecolor=color1,
	linkcolor=color1,
	bookmarksopen=false,
	pdftitle={Title},
	pdfauthor={Author},
}

\JournalInfo{1} 
\Archive{1} 

\PaperTitle{Time filtering methods for electrohydrodynamics models} 

\Authors{Li Conghui\textsuperscript{1}*} 
\affiliation{\textsuperscript{1}\textit{College of Mathematics and Statistics, Chongqing University, Chongqing, China}} 

\affiliation{*\textbf{Corresponding author}: 202106021040t@stu.cqu.edu.cn} 

\newcommand{\keywordname}{Keywords:} 

\Abstract{Electrohydrodynamics is a discipline that studies the interaction between fluid motion and electric field. Finite element method, finite difference method and other numerical simulations are effective numerical calculation methods for electrofluid dynamics models. In this paper, the finite element format of the electrofluid dynamics model is established, and the second-order convergence accuracy of the format is achieved through time filtering method. Finally, a numerical example is given to verify the convergence.}

\begin{document}

\maketitle 

\tableofcontents 

\thispagestyle{empty} 


\section{Introduction} 

 ~~~~Electrohydrodynamics (EHD) is the study of fluid mechanics in liquid dielectric under the action of electric fields, or electrodynamics in moving media, and is a cross-disciplinary study of the interaction between free charge, electric field and fluid motion. The effect of fluid motion on electric field and the effect of electric field on motion are very common phenomena, such as electroosmosis~\cite{r1}, electrophoresis ~\cite{r2} and electric field-induced instability~\cite{r3}. Electrohydrodynamic based heat transfer enhancement has the characteristics of fast response and low power consumption, and EHD has been proven to be an effective heat transfer enhancement technique in single-phase and multiphase systems. This discipline is widely used in industry, such as electro-hydrodynamic inkjet printing, to achieve ultra-high precision printing; There is also electrostatic spraying, which can be used in the industry for the spraying of automobiles, home appliances and other products, and can also be used for spraying pesticides.

In 1943, Courant solved St.Venant's torsion problem using the principle of minimum potential energy combined with piece-based continuous functions on a triangular region, which is considered to be the earliest appearance of the finite element method. With the development and application of modern computer technology, the first successful attempt of the finite element method was made in 1956. Tuner and Clough et al. successfully solved the plane stress problem with the triangular element of the finite element method when analyzing the structure of an aircraft. After continuous exploration, the basic theory and steps of finite element method gradually matured in the 1970s ~\cite{r7}. After more than 50 years of development, finite element method has gone through a series of development processes, providing new solutions to problems that could not be solved in the past, or problems that could be solved but with low accuracy. Nowadays, finite element method has been applied in various disciplines and production practices. This method to obtain numerical solutions of differential equations is not only widely used, but also can be extended for more exploration and expansion due to its characteristics~\cite{r8}, such as simple idea, clear concept and easy to understand, standardized calculation format, convenient and modular matrix expression, stable and efficient numerical calculation, etc. The value and significance of the research is very considerable. The finite element method regards the region to be solved as composed of many small adjacent subfields, which are the "element" concept in the "finite element method", and calls this process as subdivision. By means of variational method, the approximate solution on each small element is obtained in order to expect the error function of the difference between the true solution and the numerical solution to reach the minimum value and produce a stable solution.

EHD is a discipline studying fluid mechanics and electrodynamics, mainly studying the relationship between charged fluids and electric fields is a cross-disciplinary discipline. Previous studies on EHD mainly focused on single-phase flow system, and studied this problem from the aspects of stability analysis, experiment and numerical method. In recent years, with the application of EHD in complex thermal energy systems, the multiphase flow problem has gradually attracted a lot of attention. Compared with the single-phase EHD model, the multiphase EHD problem needs to consider the additional electric power and indicating charges caused by different phases at the phase interface. In physics, using some particle-based methods, dissipative particle dynamics (DPD) and lattice Boltzmann method ~\cite{r9} (LBM) are used to study the melting of EHD. In terms of numerical simulation, most research results are obtained by direct discretization of partial differential equations, such as finite difference method, finite volume method and finite element method~\cite{r7}. For the numerical simulation of complex multiphase EHD melting problem, the main difficulties are the coupling of electric field, flow field and temperature field in the model, and the change of electric field caused by the difference of physical properties of different substances.

At the end of last century, Q. Du~\cite{r12} considered the combination of ultrafine particle method and fourth-order Runge-Kutta integral. P.A. V 'Azquez et al.\cite{r13} applied the discontinuous Galerkin finite element method (DGFEM) to simulate problems related to EHD. J.L. GUERMOND et al.~\cite{r14} analyzed the convergence of a new fractional time step method for solving the incompressible Navier-Stokes equations with variable density, and established the stability of the second order variables in the form of the method. In recent years, in addition to the classical discretization method, in 2016 Kang Luo et al.~\cite{r20} established a unified lattice Boltzmann model to simulate electric convection in dielectric liquids induced by single-stage charge injection. Instead of solving complex coupled Navier-Stokes equations, charge conservation equations, and potential Poisson equations, three consistent lattice Boltzmann equations are formulated, which is a particle-based approach; Feng X et al.~\cite{r11} constructed two Gauge - Uzawa schemes of conservative scheme and convective scheme to solve the natural convection problem of variable density, and proved that the first-order versions of both schemes are unconditionally stable. MpA et al.~\cite{r10} considered several discrete schemes of variable density electrohydrodynamics models and proposed two linear time propulsion methods for the electrohydrodynamics model of fluid charge transfer in variable density media, and verified the convergence rate and energy stability of the schemes through numerical experiments.

In this paper,an electrohydrodynamics model is presented chapter 2. In the chapter 3, the finite element scheme of the model is established. Finally, the numerical solution is given in the chapter 4.


\section{EHD Model}

Define common spatial representations:
\begin{flalign*}
	& X=H_0^1(\Omega)^2=\left\{v \in H^1(\Omega)^2:v=0~~ on~~ \partial\Omega\right\}, \\
	& M=L_0^2(\Omega)=\left\{q \in L^2(\Omega):\int_{\Omega} q dx=0\right\},\\ 
	& Q=H_0^1(\Omega)=\left\{\phi \in H^1(\Omega):\phi=0~~ on~~ \partial\Omega\right\},
\end{flalign*}
here X is a two-dimensional function space, and M and Q are one-dimensional function Spaces.Let ${\mathcal{T}}_h={K_h}$ be the regular uniform triangulation of the mesh size $h$ on the calculation region $\Omega$, where $0<h<1$, the finite element space is defined as follows:
\begin{flalign*}
	&X_h=\left\{\textbf{v}_h\in X:\textbf{v}_h|_{K_h}\in (P_2(K_h))^2,\,\forall \, K_h\in {\mathcal{T}}_h\right\},\\
	&M_h=\left\{q_h \in M:q_h|_{K_h}\in P_1(K_h),\,\forall \, K_h\in {\mathcal{T}}_h\right\},\\
	&Q_h=\left\{\phi_h \in Q:\phi_h|_{K_h}\in (P_2(K_h)),\,\forall \, K_h\in {\mathcal{T}}_h\right\},
\end{flalign*}
and the finite element space $X_h \subset X$, $M_h \subset M$, $Q_h \subset Q$. Here, $P_i(K_h)$ represents the shard polynomial space with degree less than or equal to the integer $i$ on any triangulation $K_h$.

The first consideration is the electrokinetic equation 
, and the model is shown as follows:
\begin{equation}\label{mod31}
	\rho_t+\nabla\cdot(\rho\textbf{u})=0,
\end{equation}
\begin{equation}\label{mod32}
	\rho\left(\textbf{u}_t+\left(\textbf{u}\cdot\nabla\right)\textbf{u}\right)+\nabla p-\nu\Delta\textbf{u}=f,
\end{equation}
\begin{equation}\label{mod33}
	\nabla\cdot\textbf{u}=0,
\end{equation}
\begin{equation}\label{mod34}
	-\Delta\phi=\rho_e,
\end{equation}
\begin{equation}\label{mod35}
	\frac{\partial\rho_e}{\partial t}+\cdot\nabla(\rho_e\textbf{u})=\frac{1}{Pe}\Delta\rho_e-J_0\rho_e.
\end{equation}
In (\ref {mod31}) - (\ref {mod35}), $\phi $ represents potential, $\rho_e=\rho_e(\textbf {x}, t)$ represents the charge density, the fluid velocity is $\textbf{u}=\textbf{u} (\textbf {x}, t) = (u_1(\textbf{x},t), u_2(\textbf {x}, t))^T$ , $p = p (\textbf {x}, t)$ on behalf of the pressure of the fluid, $Pe$ represents the Peeler number, $J_0$ is the proportionality constant, $Re$ represents the Reynolds number, and $C_0$is the Coulomb strength constant. $\Omega\in\mathbb{R}^2$ is the convex polygon region, where it represents the incompressible fluid domain, and $T_1$ represents time and is usually positive. At the same time, give the following initial conditions on $\Omega$:
\begin{equation}
	\rho_e(\textbf{x},0)={\rho_e}_0,  
	\textbf{u}(\textbf{x},0)=\textbf{u}_0.
\end{equation}
and the boundary conditions:
\begin{equation}
	\phi|_{\Gamma}=g_1,  \rho_e|_{\Gamma}=g_2, \textbf{u}|_{\partial\Omega}=\textbf{0}.
\end{equation}
Where $\Gamma={\textbf{x}\in\partial\Omega:\textbf{g}\cdot\vec{\textbf{n}}<0}$ represents the inflow boundary, where $\vec{\textbf{n}}$ is the outer normal vector, ${\rho_e}_0, \textbf{u}_0$ both are known functions. For convenience, take $g_1=g_2=0$ and assume $\Gamma\neq\emptyset$ and $\textbf{g}\cdot\vec{\textbf{n}}=0$ on the boundary $\partial\Omega$.
Similarly, we consider the second model. Unlike the first model, where $\theta$ represents temperature, this is a temperature-dependent electrohydrodynamics model.
\begin{equation}\label{mod01}
	\frac{\partial\textbf{u}}{\partial t}+\left(\textbf{u}\cdot\nabla\right)\textbf{u}=-\nabla p+\nabla^2 \textbf{u}+\left(\frac{T}{M}\right)^2qC(-\nabla\phi),
\end{equation}
\begin{equation}\label{mod02}
	\nabla\cdot\left(\frac{1}{C}\nabla\phi\right)=-q,
\end{equation}
\begin{equation}\label{mod03}
	\frac{\partial q}{\partial t}+\nabla\cdot\left[\left(\frac{T}{M^2}(-\nabla\phi)+\textbf{u}\right)q\right]=\alpha\nabla^2 q,
\end{equation}
\begin{equation}\label{mod04}
	\frac{\partial \theta}{\partial t}+\textbf{u}\cdot\nabla\theta=\frac{1}{Pr}\nabla^2\theta-\frac{1}{St}\frac{\partial f_l}{\partial t},
\end{equation}
\begin{equation}\label{mod05}
	\nabla\cdot\textbf{u}=0.
\end{equation}

\section{Finite element methods for the model}
First, the first model is analyzed and its variational form is obtained~\cite{r21}. Using (\ref {mod31}) - (\ref {mod35}) belong to different space unknown functions in the corresponding test function, for $\forall (v,q,\omega,\varphi,\psi)\in(X,M,Q,Q,Q)$, solve $(u,p,\rho,\rho_e,\phi)\in (X,M,Q,Q)$, such that for $\forall t\in(0,T]$, there is

Step 1.~~Solve $\tilde{\rho}_h^{n+1}$
\begin{align}
	&\left(\frac{\tilde{\rho}_h^{n+1}-\rho_h^{n}}{\Delta t},w_h\right)+\left(\nabla\cdot(\tilde{\rho}_h^{n+1}\tilde{\textbf{u}}_h^{n+1}),w_h\right)\notag\\
	&-\frac{1}{2}\left(\tilde{\rho}_h^{n+1}\nabla\cdot\tilde{\textbf{u}}_h^{n+1},w_h\right)=0,
\end{align}

Step 2.~~Solve $\tilde{\textbf{u}}_h^{n+1}$,$p_h^{n+1}$
\begin{align}
		&\left(\rho_h^{\star}\frac{\tilde{\textbf{u}}_h^{n+1}-\textbf{u}_h^{n}}{\Delta t},\textbf{v}_h\right)
		+b_1\left(\tilde{\rho}_h^{n+1}u_h^{\star},\tilde{\textbf{u}}_h^{n+1},\textbf{v}_h\right)\notag\\
		&-\left(\tilde{p}_h^{n+1},\nabla\cdot\textbf{v}_h\right)
		+\left(\nabla\cdot\tilde{\textbf{u}}_h^{n+1},q_h\right)\notag\\
		&+\nu\left(\nabla\tilde{\textbf{u}}_h^{n+1},\nabla\textbf{v}_h\right)\notag\\
		&+\frac{1}{4}b_2\left(\tilde{\rho}_h^{n+1}\tilde{\textbf{u}}_h^{n+1},\tilde{\textbf{u}}_h^{n+1},\textbf{v}_h\right)=(f^{n+1},\textbf{v}_h),
\end{align}
\begin{align}
	&\left(\frac{\tilde{\rho_e}_h^{n+1}-{\rho_e}_h^{n}}{\Delta  t},\xi_h\right)+\left(\nabla\cdot(\tilde{{\rho_e}}_h^{n+1}\textbf{u}_h^{\star}),\xi_h\right)\notag\\
	&-\frac{1}{Pe}\left(\nabla\tilde{\rho_e}_h^{n+1},\nabla\xi_h\right)+J_0\left(\tilde{\rho_e}_h^{n+1},\xi_h\right)=0,
\end{align}

Step 3.~~Solve $\tilde{\phi_h^{n+1}}$
\begin{equation}
	\left(\nabla\tilde{\phi}_h^{n+1},\nabla\psi_h\right)=(\tilde{\rho_e}_h^{n+1},\psi_h).
\end{equation}

Step 4.~~Time filtering
\begin{equation}
	{\rho}_h^{n+1}=\tilde{\rho}_h^{n+1}-\frac{1}{3}\left(\tilde{\rho}_h^{n+1}-2{\rho_e}_h^{n}+{\rho}_h^{n-1}\right),
\end{equation}
\begin{equation}
	\textbf{u}_h^{n+1}=\tilde{\textbf{u}}_h^{n+1}-\frac{1}{3}\left(\tilde{\textbf{u}}_h^{n+1}-2\textbf{u}_h^{n}+\textbf{u}_h^{n-1}\right),
\end{equation}
\begin{equation}
	{\rho_e}_h^{n+1}=\tilde{\rho_e}_h^{n+1}-\frac{1}{3}\left(\tilde{\rho_e}_h^{n+1}-2{\rho_e}_h^{n}+{\rho_e}_h^{n-1}\right),
\end{equation}
\begin{equation}
	\phi_h^{n+1}=\tilde{\phi}_h^{n+1}-\frac{1}{3}\left(\tilde{\phi}_h^{n+1}-2\phi_h^{n}+\phi_h^{n-1}\right),
\end{equation}
here,
$\rho_h^ {\star}=2\rho_h^{n}-\rho_h^{n-1}$, $\textbf{u}_h^{\star}=2\textbf{u}_h^{n}-\textbf{u}_h^{n-1}$, ${\rho_e}_h^{\star}=2{\rho_e}_h^{n}-{\rho_e}_h^{n-1}$.

Then the second model is analyzed and its variational form is obtained. Solve$(\textbf{u},p,q,\phi,\theta)$ such that $\forall t\in(0,T]$,

\begin{align}
	&\left(\frac{\partial \textbf{u}}{\partial t},\textbf{v}\right)+\left((\textbf{u}\cdot\nabla)\textbf{u},\textbf{v}\right)+(\nabla p,\textbf{v})\notag\\
	&+(\nabla\textbf{u},\nabla\textbf{v})+\left(\frac{T}{M}\right)^2C\left(q\nabla\phi,\textbf{v}\right)=0,
\end{align}
\begin{equation}
	\frac{1}{C}\left(\nabla\phi,\nabla\psi\right)-(q,\psi)=0,
\end{equation}
\begin{align}
	&\left(\frac{\partial q}{\partial t},\xi\right)-\frac{T}{M^2}\left(\left(-\nabla\phi\right)q,\nabla\xi\right)-\left(\textbf{u}q,\nabla\xi\right)\notag\\
	&+\alpha(\nabla q,\nabla \xi)=0,
\end{align}
\begin{equation}
	\left(\frac{\partial \theta}{\partial t},t\right)+(\textbf{u}\cdot\nabla\theta,t)+\frac{1}{Pr}\left(\nabla\theta,\nabla t\right)=0,
\end{equation}
\begin{equation}
	\left(\nabla\cdot\textbf{u},w\right)=0,
\end{equation}

here,$\textbf{v},\psi,\xi,t,w$ are test functions.

Then the finite element format of the second model is established. Given $\textbf{u}_h^{n}=\textbf{u}^{0}$, $p_h^n=p^0$, $q_h^{n}=q^{0}$, $\phi_h^n=\phi^0$, $\theta_h^{n}=\theta^{0}$, $n=0,1,..., N-1$,solve ($\textbf{u}_h^{n+1}$, $p_h^{n+1}$, $q_h^{n+1}$, $\phi_h^{n+1}$, $\theta_h^{n+1}$),such that

Step 1. solve $u_h^{n+1}$and$p_h^{n+1}$:
\begin{align}
	&\left(\frac{\textbf{u}_h^{n+1}-\textbf{u}_h^{n}}{\Delta t},\textbf{v}_h\right) +\left(\left(\textbf{u}_h^{n}\cdot\nabla\right)\textbf{u}_h^{n+1},\textbf{v}_h\right)\notag\\
	&+(\nabla p_h^{n+1},\textbf{v}_h)
	+(\nabla\textbf{u}_h^{n+1},\nabla\textbf{v}_h)\notag\\
	&+\left(\frac{T}{M}\right)^2 C(q_h^{n}\nabla\phi_h^{n},\textbf{v}_h)=0,
\end{align}
\begin{equation}
	\left(\nabla\cdot\textbf{u}_h^{n+1},w_h\right)=0,
\end{equation}
	
	Step 2. solve $q_h^{n+1}$and$\phi_h^{n+1}$:
	\begin{align}
		&\left(\frac{q_h^{n+1}-q_h^{n}}{\Delta t},\xi_h\right)-\frac{T}{M^2}\left(\left(-\nabla\phi_h^{n+1}\right)q_h^{n},\nabla\xi_h\right)\notag\\
		&-\left(\textbf{u}_h^{n}q_h^n,\nabla\xi_h\right)+\alpha(\nabla q_h^{n+1},\nabla \xi_h)=0,
	\end{align}
	
	\begin{equation}
		(q_h^{n+1},\psi_h)-\frac{1}{C}\left(\nabla\phi_h^{n+1},\nabla\psi\right)=0,
	\end{equation}
	
	Step 3. solve $\theta_h^{n+1}$:
	\begin{align}
		&\left(\frac{\theta_h^{n+1}-\theta_h^{n}}{\Delta t},t_h\right)+\left(\textbf{u}_h^n\cdot\nabla\theta_h^{n+1},t_h\right)\notag\\
		&-\frac{1}{Pr}\left(\nabla\theta_h^{n+1},\nabla t_h\right)=0.
	\end{align}


\section{Numerical experiment}
Considering the known example of true solution, regional $\Omega$ to [0, 1] x [0, 1], equations (\ref {mod31}) - (\ref {mod35}) take the functions as true solutions
\begin{equation*}
	\rho=2+x cos(sin(t))+y sin(sin(t)),
\end{equation*}
\begin{equation*}
	\textbf{u}
	=
	\left(
	\begin{matrix}
		u_1
		\\u_2
	\end{matrix}
	\right)
	=
	\left(
	\begin{matrix}
		-y cos(t)
		\\x cos(t)
	\end{matrix}
	\right),
\end{equation*}
\begin{equation*}
	p=sin(x) sin(y) sin(t),
\end{equation*}
\begin{equation*}
	\rho_e=2 sin(x) sin(y) sin(t),
\end{equation*}
\begin{equation*}
	\phi=sin(x) sin(y) sin(t).
\end{equation*}
Furthermore, we implemented the codes using the
software package FreeFEM++.The results of each variable in the numerical experiment are shown in the table \ref{tab1}.
\begin{table*}[htb]
	\caption[1]{Error and convergence order of true solution for variable density model}
	\centering
	\renewcommand{\arraystretch}{1.5}
	\setlength{\tabcolsep}{1.4pt}
	\begin{tabular}{ccccccccccc}
		\toprule[1.5pt]
		N & $\Vert e_{\rho}\Vert_0$ & Order &$\Vert e_{\textbf{u}}\Vert_0$ & Order & $\Vert e_p\Vert_0$ & Order & $\Vert e_{\rho_e}\Vert_0$ & Order & $\Vert e_{\phi}\Vert_0$ & Order \\
		\midrule[1pt]
		4 & 1.684e-02 & - & 2.424e-02 & - & 9.251e-02 & - & 1.008e-02 & - & 5.730e-03 \\ 
		8 & 6.317e-03 & 1.4146 & 7.115e-03 & 1.7682 & 1.629e-02 & 2.5058 & 3.015e-03 & 2.5058 & 1.624e-03 & 2.5058 \\ 
		16 & 1.874e-03 & 1.7532 & 1.825e-03 & 1.9628 & 4.677e-03 & 1.8001 & 7.962e-04 & 1.8001 & 4.241e-04 & 1.8001 \\
		32 & 5.037e-04 & 1.8955 & 4.523e-04 & 2.0129 & 1.190e-03 & 1.9702 & 2.041e-04 & 1.9702 & 1.078e-04 & 1.9702 \\
		64 & 1.302e-04 & 1.9513 & 1.112e-04 & 2.0245 & 3.165e-04 & 1.9111 & 5.127e-05 & 1.9513 & 2.753e-05 & 1.9696 \\
		128 & 3.310e-05 & 1.9762 & 2.736e-05 & 2.0225 & 9.930e-05 & 1.6725 & 1.289e-05 & 1.9762 & 6.918e-06 & 1.9926 \\
		\bottomrule[1.5pt]
	\end{tabular}
	\label{tab1}
\end{table*}

For the second temperature-dependent model, the true solutions are as follows:
\begin{equation}
	\textbf{u}
	=
	\left(
	\begin{matrix}
		u_1
		\\u_2
	\end{matrix}
	\right)
	=
	\left(
	\begin{matrix}
		-y cos(t)
		\\x cos(t)
	\end{matrix}
	\right),
\end{equation}
\begin{equation}
	p=sin(x) sin(y) sin(t),
\end{equation}
\begin{equation}
	q=2 sin(x) sin(y) sin(t),
\end{equation}
\begin{equation}
	\phi=sin(x) sin(y) sin(t),
\end{equation}
\begin{equation}
	\theta=-y cos(t)+x cos(t).
\end{equation}
The results in Table \ref{tab2} can be obtained by writing programs.
\begin{table*}[htb]
	\centering
	\caption{Numerical results of a temperature dependent model}
	\renewcommand{\arraystretch}{1.5}
	\setlength{\tabcolsep}{1.8pt}
	\begin{tabular}{ccccccccccc}
		\toprule[1.5pt]
		N & $\Vert e_{\textbf{u}}\Vert_0$ & Order & $\Vert e_p\Vert_0$ & Order & $\Vert e_q\Vert_0$ & Order & $\Vert e_{\phi}\Vert_0$ & Order & $\Vert e_{\theta}\Vert_0$ &Order\\
		\midrule[1pt]
		4 & 7.2199e-04 & - & 4.6444e-02 & - & 6.5729e-03 & - & 3.2538e-04 & - & 4.2763e-03 \\ 
		8 & 3.7046e-04 & 0.9627 & 2.4509e-02 & 0.9222 & 3.3066e-03 & 0.9912 & 1.6397e-04 & 0.9887 & 2.3874e-03 & 0.8409 \\ 
		16 & 1.8730e-04 & 0.9840 & 1.2550e-02 & 0.9656 & 1.6667e-03 & 0.9884 & 8.2417e-05 & 0.9924 & 1.2414e-03 & 0.9435 \\
		32 & 9.4060e-05 & 0.9937 & 6.3374e-03 & 0.9857 & 8.3551e-04 & 0.9962 & 4.1265e-05 & 0.9980 & 6.3366e-04 & 0.9702 \\
		64 & 4.7142e-05 & 0.9966 & 3.1830e-03 & 0.9935 & 4.1815e-04 & 0.9986 & 2.0653e-05 & 0.9985 & 3.1988e-04 & 0.9862 \\
		128 & 2.3607e-05 & 0.9978 & 1.5958e-03 & 0.9961 & 2.0919e-04 & 0.9993 & 1.0331e-05 & 0.9994 & 1.6071e-04 & 0.9931 \\
		\bottomrule[1.5pt]
	\end{tabular}
	\label{tab2}
\end{table*}


\phantomsection
\section{Conclusions} 

In this paper, the finite element formats of temperature-independent and temperature-dependent current EHD models are established respectively, and their convergence order is verified by numerical experiments.

\newpage
\phantomsection
\bibliographystyle{unsrt}
\bibliography{sample.bib}


\end{document}